\documentclass[twoside,12pt]{article}
\usepackage{amssymb}
\usepackage{amsmath}
\usepackage[dvips]{graphicx}

\newtheorem{theorem}{Theorem}

\newtheorem{corollary}{Corollary}
\newtheorem{remark}{Remark}
\newtheorem{proposition}{Proposition}

\title{Optimal temperature distribution for a nonisothermal Cahn--Hilliard system in two dimensions with source term and double 
obstacle potential}

\def\enne{{\mathbb{N}}}
\def\erre{{\mathbb{R}}}
\def\nn{{\bf n}}
\def\dn{{\partial_{\bf n}}}
\def\dt{{\partial_t}}
\def\phi{\varphi}
\def\dtt{{\partial_{tt}}}
\def\indi{{I_{[-1,1]}}}
\def\sindi{{\partial \indi}}

\def\kuno{{\kappa_1}}
\def\kdue{{\kappa_2}}

\def\J{{\cal J}}
\def\S{{\cal S}}

\def\Uad{{\cal U}_{\rm ad}}

\def\calU{{\cal U}}
\def\non{{\nonumber}}
\def\Ldue{{L^2(\Omega)}}
\def\Huno{{H^1(\Omega)}}
\def\Hdue{{H^2(\Omega)}}

\def\iot {\int_0^t}

\def\intQt{\int_{Q_t}}

\def\iO{\int_\Omega}

\def\Liq{{L^\infty(Q)}}

\def\hal{h_\alpha}
\def\dhal{h_\alpha'}
\def\phial{\varphi_\alpha}
\def\mual{\mu_\alpha}
\def\wal{w_\alpha}
\def\phialn{\varphi_{\alpha_n}}
\def\mualn{\mu_{\alpha_n}}
\def\xialn{h_{\alpha_n}'(\phialn)}
\def\waln{w_{\alpha_n}}
\def\pal{p_\alpha}
\def\qal{q_\alpha}
\def\ral{r_\alpha}
\def\Sal{{\cal S}_\alpha}
\def\Saln{{\cal S}_{\alpha_n}}

\def\vbar{\overline v}

\def\psibar{\overline\psi}

\def\eps{{\varepsilon}}

\def\qed{{\hfill$\square$}}
\def\calN{{\cal N}}
\def\Vp{{V^*}}
\def\vp{{\varphi}}

\date{~}
\begin{document}

\maketitle
\begin{center}
\vskip-1.5cm
{\large\sc Pierluigi Colli$^{(1)}$}\\
{\normalsize e-mail: {\tt pierluigi.colli@unipv.it}}\\[0.25cm]
{\large\sc Gianni Gilardi $^{(1)}$}\\
{\normalsize e-mail: {\tt gianni.gilardi@unipv.it}}\\[0.25cm]
{\large\sc Andrea Signori$^{(2)}$}\\
{\normalsize e-mail: {\tt andrea.signori@polimi.it}}\\[0.25cm]
{\large\sc J\"urgen Sprekels$^{(3)}$}\\
{\normalsize e-mail: {\tt juergen.sprekels@wias-berlin.de}}\\[.5cm]
$^{(1)}$
{\small Dipartimento di Matematica ``F. Casorati'', Universit\`a di Pavia}\\
{\small and Research Associate at the IMATI -- C.N.R. Pavia}\\
{\small via Ferrata 5, I-27100 Pavia, Italy}\\[.3cm] 
$^{(2)}$
{\small Dipartimento di Matematica, Politecnico di Milano}\\
{\small via E. Bonardi 9, I-20133 Milano, Italy}
\\[.3cm] 
$^{(3)}$
{\small Department of Mathematics}\\
{\small Humboldt-Universit\"at zu Berlin}\\
{\small Unter den Linden 6, D-10099 Berlin, Germany}\\
{\small and}\\
{\small Weierstrass Institute for Applied Analysis and Stochastics}\\
{\small Mohrenstrasse 39, D-10117 Berlin, Germany}\\[10mm]
\end{center}
\date{}

\pagestyle{myheadings}
\newcommand\testopari{\sc Colli -- Gilardi -- Signori -- Sprekels}
\newcommand\testodispari{\sc Optimal Control of a nonisothermal Cahn--Hilliard type model}
\markboth{\testopari}{\testodispari}

{\small
\begin{abstract}
\noindent In this note, we study the optimal control of a  nonisothermal phase field system of Cahn--Hilliard type that constitutes an 
extension of the classical Caginalp model for nonisothermal phase transitions with a conserved order parameter. It couples a Cahn--Hilliard type equation with source term for the order parameter with the universal balance law of internal energy. In place of the standard 
Fourier form, the constitutive law of the heat flux is assumed in the form given by the theory developed by Green and Naghdi, which 
accounts for a possible thermal memory of the evolution. This has the consequence that the balance law of internal energy becomes
a second-order in time equation for the {\it thermal displacement} or {\it freezing index}, that is, a primitive with respect to time of
the temperature. Another particular feature of our system is the presence of the source term in the equation for the order parameter, which
entails further mathematical difficulties because the mass conservation of the order parameter is no longer satisfied. 
In this paper, we study the case that the double-well potential driving the evolution of the phase transition is given by the nondifferentiable
double obstacle potential, thereby complementing recent results obtained for the differentiable cases of regular and logarithmic potentials.
Besides existence results, we derive first-order necessary optimality conditions for the control problem.
The analysis is carried out by employing the so-called {\it deep quench approximation} in which the nondifferentiable double obstacle potential is
approximated by a family of potentials of logarithmic structure for which meaningful first-order necessary optimality conditions in terms of 
suitable adjoint systems and variational inequalities are available. Since the results for the logarithmic potentials crucially depend 
on the validity of the so-called {\it strict separation property} which is only available in the spatially two-dimensional situation,
our whole analysis is restricted to the two-dimensional case.  

\vskip3mm



\end{abstract}
}

\section{Introduction}
\label{SEC:INTRO}
\setcounter{equation}{0}

\noindent 
Let $\Omega\subset\erre^2$ be some open, bounded, and connected set 
having a smooth boundary $\Gamma:=\partial\Omega$ and the outward unit normal field~$\,\nn$. 
Denoting by $\dn$ the directional derivative in the direction of $\nn$, and putting, with a fixed final time $T>0$,
$$
  Q_t:=\Omega\times(0,t) \,\mbox{ and }\,\Sigma_t:=\Gamma\times(0,t)\,\mbox{ for $\,t\in(0,T] $, \quad $Q:=Q_T$,\quad
  $\Sigma:=\Sigma_T $},
$$ 
we study in this paper as {\it state system} the following initial-boundary value problem:
\begin{align}
  & \dt\phi - \Delta\mu + \gamma \phi
  = f  
  \quad \mbox{in $Q$},
  \label{Iprima}
  \\
  & \mu = - \Delta\phi + \xi+F'(\phi) + a - b \dt w,\quad\xi\in\sindi(\phi),
  \quad \mbox{in $Q$},
  \label{Iseconda}
  \\
  & \dtt w - \Delta(\kuno \dt w + \kdue w) + \lambda \dt\phi
  = u
  \quad\mbox{in $Q$},
  \label{Iterza}
	\\
  &\dn\phi
  = \dn\mu
  = \dn(\kuno \dt w + \kdue w)
  = 0 
  \quad \mbox{on $\Sigma$},
  \label{Ibc}
  \\
  &\phi(0) = \phi_0, 
  \quad 
  w(0) = w_0,
  \quad 
  \dt w(0) = w_1,
  \quad \mbox{in $\Omega.$}
  \label{Icauchy}
\end{align}
The {\it cost functional} under consideration is given by
\begin{align}
	&\J ((\phi, w), u)
	:= {}
		\frac {\beta_1}2 \int_Q |\phi - \phi_Q|^2
	+ \frac {\beta_2}2  \iO|\phi(T) - \phi_\Omega|^2
	\non
	\\ 
	& {}
	+ \frac {\beta_3}2 \int_Q |w - w_Q|^2
	+ \frac {\beta_4}2 \iO |w(T) -  w_\Omega|^2
	\non
	\\ 
	& {}
	+ \frac {\beta_5}2 \int_Q |\dt w -  w'_Q|^2
	+ \frac {\beta_6}2 \iO |\dt w(T) -  w'_\Omega|^2
	+ \frac{\nu}{2} \int_Q |u|^2,
	\label{cost}
\end{align}
with nonnegative constants $\beta_i$, $1\le i\le 6$, and $\nu$, which are not all zero, and where $\phi_\Omega, w_\Omega, w'_\Omega\in\Ldue$ and
$\phi_Q,w_Q,w'_Q\in L^2(Q)$ denote given target functions.

For the distributed control variable $u$, we choose as control space 
\begin{equation}
  \calU := L^\infty(Q),
  \label{defU}
\end{equation}
and the related set of admissible controls
is given by
\begin{equation}
	\label{Uad}
	\Uad : = \big\{ u \in {\cal U}: u_{\rm min} \leq u \leq u_{\rm max} \,\mbox{ a.e. in }\,Q \big\},
\end{equation}
where we generally assume throughout the paper that 
\begin{equation}
\label{uminmax}
 u_{\rm min},u_{\rm max}\in\Liq \quad\mbox{and}\quad u_{\rm min}\le u_{\rm max}\,\mbox{ a.e. in }\,Q.
\end{equation}
 In particular,
$\Uad$ is bounded in $\Liq$.

In summary, the control problem under investigation can be reformulated as follows:

\vspace{2mm}
\noindent    
{\bf (P)} \,\,\,$ \min_{u \in \Uad}\J ((\phi, w), u)$ \,\,subject to the constraint that $(\phi,\mu, \xi,w)$ 
 solves  
 the state \linebreak \hspace*{9.5mm} system (\ref{Iprima})--(\ref{Icauchy}).

\vspace{2mm}
Let us now spend some comments on the state system (\ref{Iprima})--(\ref{Icauchy}), which is a formal extension of the nonisothermal Cahn--Hilliard system introduced by Caginalp in \cite{Cag2}
to model the phenomenon of nonisothermal phase segregation in binary mixtures (see also \cite{Cag1, Cag3} and the derivation in \cite[Ex.~4.4.2, 
(4.44), (4.46)]{BS}); it corresponds to the Allen--Cahn counterpart analyzed 
in \cite{CSS3}. The unknowns in the state system have the following physical meaning: $\phi$ is a normalized difference between 
the volume fractions
of pure phases in the binary mixture (the dimensionless {\em order parameter} of the phase transformation, which should attain its values in 
the physical interval $[-1,1]$, where the extremes represent the pure phases of the mixture), 
$\mu$ is the associated 
{\em chemical potential}, and $\,w\,$ is the so-called {\em thermal displacement} (or 
{\em freezing index}), which is directly connected to the temperature $\vartheta$ (which in the case of the Caginalp model is actually 
a temperature difference) through the relation
\begin{equation}
	\label{thermal_disp}	
	w (\cdot , t)  = w_0 + \iot \vartheta(\cdot, s) \,ds, \quad t \in[0,T].
\end{equation}
Moreover, $\kuno$ and $\kdue$ in (\ref{Iterza}) stand for prescribed positive coefficients related to the heat flux, which 
is here assumed in the Green--Naghdi form
(see~\cite{GN91,GN92,GN93,PG09})
\begin{equation}\label{flux}
\mathbf q=-\kappa_1 \nabla (\dt w )- \kappa_2 \nabla w,
\end{equation} 
which accounts for a possible previous
thermal history of the phenomenon. Moreover,
$\gamma$ is a positive physical constant related to the intensity of the mass absorption/production of the source,
where the source term in (\ref{Iprima}) is $\,S:=f- \gamma \varphi$. This term reflects the fact that the system may not be isolated and a loss or production of mass is possible, which happens, e.g., in numerous liquid-liquid phase segregation problems that arise in cell biology~\cite{Bio}
and in tumor growth models~\cite{GLSS}. Notice that the presence of the source term entails that the property of mass conservation of the order
parameter is no longer valid; in fact, from (\ref{Iprima}) it directly follows that the mass balance has the form
\begin{equation}
\label{masscons} 
\frac d{dt} \, \Bigl(\frac 1 {|\Omega|}\iO \phi(t) \Bigr)
	= 
	\frac 1 {|\Omega|} \iO {S(t)}, \quad\mbox{for a.e. $t\in(0,T)$},
\end{equation}
where $\,|\Omega|\,$ denotes the Lebesgue measure of $\Omega$.

In addition to the quantities already introduced,
$\lambda$ stands for the latent heat of the phase transformation, $a, b$ are physical constants, and the control variable
$\,u\,$ is  a distributed heat source/sink. Besides,  $\phi_0,w_0,$ and $w_1$ indicate some given initial values.  
Moreover, the function $F$, whose derivative appears in (\ref{Iseconda}), is assumed to be concave, typically of the form $\,F(r)=c_1-c_2r^2$
 with $c_1\in\erre$, $c_2>0$,  while $\,\sindi\,$
denotes the subdifferential of the indicator function $\,\indi\,$ of the real interval $[-1,1]$, which is given by 
\begin{equation}
\label{indi} 
\indi(r)=0\quad\mbox{if $|r|\le 1$}\quad\mbox{and} \quad\indi(r)=+\infty\quad\mbox{if $|r|>1$}\,.
\end{equation}
The potential
\begin{equation}
\label{F2obs}
F_{2obs}(r)=\indi(r)+F(r),
\end{equation}
with $\,F\,$ given as above, is then the typical double obstacle potential.

The state system (\ref{Iprima})--(\ref{Icauchy}) was recently analyzed in \cite{CGSS3} concerning well-posedness and
regularity (see the results cited below in Section 2); in
 \cite{CGSS4} the corresponding optimal control problem {\bf (P)} has been analyzed for the simpler differentiable case
when the indicator function $\indi$ occurring in (\ref{Iseconda}) is replaced by either a regular function or by a logarithmic expression of the form
\begin{equation}
  \label{defhal}
  h_{\alpha}(r) := \alpha\,h(r),
\end{equation}
with $\alpha>0$, where
\begin{equation}
\label{defh}
h(r)=	
  \left\{
    \begin{array}{ll}
      (1+r)\ln (1+r)+(1-r)\ln (1-r)  
      & \quad \hbox{if $|r|<1$}
      \\
			2\,\ln(2)
			&\quad \hbox{if $r\in\{-1,1\}$}\\
      +\infty
      & \quad \hbox{if $|r|>1$}
	\end{array}
  \right.\,.
  \end{equation}
Clearly, in this case the subdifferential inclusion (\ref{Iseconda}) has to be replaced by the equation
\begin{equation}
\label{Isecondanew} 
\mu=-\Delta\varphi+h_\alpha'(\varphi)+F'(\varphi)+a-b \dt w\,.
 \end{equation}
For such logarithmic nonlinearities, in \cite{CGSS4} results concerning existence of optimal controls, Fr\'echet differentiability of
the control-to-state operator, and meaningful first-order necessary optimality conditions  
(in terms of the associated adjoint state problem and
variational inequality) have been derived, at least in the spatially two-dimensional situation. In this paper, we complement the
results of \cite{CGSS4} by investigating the optimal control problem for the nondifferentiable double obstacle case. 
While the existence of optimal controls is not too difficult to show, the derivation of 
first-order necessary optimality conditions is a much more challenging task, since the existence of appropriate Lagrange multipliers
cannot be derived from the standard theory. We therefore 
employ the so-called {\it deep quench approximation}, which has been successfully applied in a number of Allen--Cahn or 
Cahn--Hilliard systems (see, e.g., \cite{CFS,CSS4,CGS1,CGS2,CGS3,CGS4,CGS26,Sig1}). The general strategy of this approach is the following.    
At first, we observe the following facts: it is readily seen that
\begin{equation}
\label{halcon}
\lim_{\alpha\searrow0}\,\hal(r)=\indi(r)\quad\forall\,r\in\erre .
\end{equation}
Moreover, $h'(r)=\ln(\frac{1+r}{1-r})$\, and \,$h''(r)=\frac 2{1-r^2}$, and thus
\begin{align}
\label{halscon}
&\lim_{\alpha\searrow0} \dhal(r)=0 \quad\,\hbox{for all }\,r\in (-1,1),
\nonumber\\
& \lim_{\alpha\searrow0}\Bigl(\lim_{r\searrow -1}\,\dhal (r)\Bigr)=-\infty, \quad
\lim_{\alpha\searrow0}\Bigl(\lim_{r\nearrow 1}\,\dhal (r)\Bigr)=+\infty.
\end{align}
Hence, we may regard the graphs of the single-valued $\alpha$-dependent functions $\dhal$ over 
the interval $(-1,1)$ as approximations to the graph of the subdifferential~$\partial\indi$ from the interior of $(-1,1)$ (in contrast to the
exterior approximation obtained via
the Moreau--Yosida approach). 

In view of the convergence properties (\ref{halcon}) and (\ref{halscon}), it is to be expected that the solutions to the approximating system 
(\ref{Iprima}), (\ref{Isecondanew}), (\ref{Iterza})--(\ref{Icauchy}) converge in a suitable topology to the solution of the state system
(\ref{Iprima})--(\ref{Icauchy}) as $\alpha\searrow0$, and a similar behavior ought to be true for the corresponding minimizers of the 
associated optimal control problems. It is then hoped that it is possible to pass to the limit as $\alpha\searrow0$ in the first-order
necessary optimality conditions for the approximating control problems in order to establish first-order conditions also for the double
obstacle case. It turns out that this general strategy works with suitable modifications.
Let us stress at this point that our approach makes use of the results obtained for the logarithmic case investigated in \cite{CGSS4};    
since in that case the derivation of differentiability 
properties of the associated control-to-state operator was only possible under the premise that the order parameter $\,\phi\,$ satisfies
the so-called \emph{strict separation property} (meaning that $\,\phi\,$ attains its values in a compact subset of $(-1,1)$), and since
this property could only be shown in the spatially two-dimensional case, our analysis does not apply to three-dimensional domains
$\,\Omega$.

The plan of the paper is as follows.
The next section is devoted to collect previous results concerning the well-posedness of the state system.
Then, in Section 3 and Section 4, we investigate the convergence properties of the deep quench approximations and of the associated optimal controls.
The final section brings the derivation of first-order necessary conditions of optimality for the problem {\bf (P)}
by employing the strategy explained above.

Prior to this, let us fix some notation.
For any Banach space $X$, we denote by \,$\|\,\cdot\,\|_X$, $X^*$, and $\langle\,\cdot\,,\, \cdot\,\rangle_X$,  
the corresponding norm, its dual space, and  the related duality pairing between $X^*$ and~$X$.
For two Banach spaces $X$ and $Y$ that are both continuously embedded in some topological vector space~$Z$, we introduce the linear space
$X\cap Y$, which becomes a Banach space when equipped with its natural norm 
$\|v\|_{X\cap Y}:=\|v\|_X+\|v\|_Y\,$, for $v\in X\cap Y$. 
A special notation is used for the standard Lebesgue and Sobolev spaces defined on $\Omega$.
For every $1 \leq p \leq \infty$ and $k \geq 0$, they are denoted by $L^p(\Omega)$ and $W^{k,p}(\Omega)$, 
with the associated norms $\|\,\cdot\,\|_{L^p(\Omega)}=\|\,\cdot\,\|_{p}$ and $\,\|\,\cdot\,\|_{W^{k,p}(\Omega)}$, respectively. 
If $p=2$, they become Hilbert spaces, and we employ the standard convention $H^k(\Omega):= W^{k,2}(\Omega)$. 
For convenience, we also~set
$$
  H := \Ldue , \quad  
  V := \Huno,   \quad
  W := \{v\in\Hdue: \ \dn v=0 \,\mbox{ on $\,\Gamma$}\}.
$$
For simplicity, we use the symbol $\|\,\cdot\,\|$ for the norm in $H$ and in any power thereof,
and we denote by $(\,\cdot\,,\,\cdot\,)$ and $\langle\,\cdot\,,\,\cdot\rangle$ the inner product in $H$ and the dual
pairing between $V^*$ and $V$.
Observe that the embeddings 
 $\, W \subset V \subset H \subset V^*\,$
are dense and compact. 
As usual, $H$ is identified with a subspace of $V^*$ to have the Hilbert triplet $(V,H,V^*)$ along with the identity
$$
	\langle u,v \rangle =(u,v)	\quad\mbox{for every $u\in H$ and $v\in V$.}
$$

Next, for a generic element $v\in V^*$,  we define its generalized mean value $\vbar$ by 
\begin{equation}
  \vbar := \frac 1{|\Omega|} \, \langle v , {\bf 1} \rangle,
  \label{defmean}
\end{equation}
where ${\bf 1}$ stands for the constant function that takes the value $1$ in $\Omega$.
It is clear that $\vbar$ reduces to the usual mean value if $v\in H$.
The same notation $\vbar$ is employed also if $v$ is a time-dependent function.

To conclude, for normed spaces $\,X\,$ and $\,v\in L^1(0,T;X)$, we define the convolution products
\begin{equation}
	({\bf 1} * v)(t):=\int_0^t v(s)\,{\rm d}s,
%	\label{conv}
	\quad 
	({\bf 1} \circledast v)(t):=\int_t^T v(s)\,{\rm d}s ,
	\qquad \hbox{$t \in[0,T]$.}
	\label{conv:back}	
\end{equation}

\section{General assumptions and the state system}

For the remainder of this paper, we make the following general assumptions besides (\ref{Uad}) and (\ref{uminmax}).
\begin{description}
\item[(A1)] \,\,The structural constants $\gamma$, $a$, $b$, $\kuno$, $\kdue$, and $\lambda$ are positive.
\item[(A2)] \,\,It holds $F\in C^3(\erre)$, and $F'$ is Lipschitz continuous on $\erre$.
\item[(A3)] \,\,$f\in H^1(0,T;H)\cap L^\infty(Q)$, \,\,\,$w_0\in V$, \,\,\,$w_1\in W$.
\item[(A4)] \,\,$\phi_0\in H^4(\Omega)\cap W$ satisfies  $\,\Delta\phi_0\in W$, and, with $\rho:=\frac 1\gamma\,\|f\|_{L^\infty(Q)}$, 
we assume that 
all of the quantities 
$$
 \inf_{x\in\Omega}\phi_0(x) , \ 
  \sup_{x\in\Omega}\phi_0(x), \
  - \rho - (\overline{\phi_0})^- \,, \ \rho + (\overline{\phi_0})^+
$$
belong to the interior of $(-1,1)$,
where $(\cdot)^+$  and $(\cdot)^-$ denote the positive and negative part functions, respectively. 
\end{description}

The analysis of the systems (\ref{Iprima})--(\ref{Icauchy}) and (\ref{Iprima}), (\ref{Isecondanew}), (\ref{Iterza})--(\ref{Icauchy}) 
has been the subject of investigation in \cite{CGSS3}.
As a special case of \cite[Thm.~2.2]{CGSS3}, we have the following result for the initial-boundary value problem (\ref{Iprima})--(\ref{Icauchy}).
\begin{theorem}
\label{SS0}
Suppose that {\rm (\ref{Uad})}, {\rm (\ref{uminmax})} and {\bf (A1)}--{\bf (A4)} are fulfilled. Then the state system 
{\rm (\ref{Iprima})--(\ref{Icauchy})} has for every $u\in\Uad$ 
a weak solution $(\phi,\mu,\xi,w)$ in the 
following sense: it holds
\begin{align}
\label{regphi0}
&\phi\in H^1(0,T;V)\cap L^\infty(0,T;W)\cap L^\infty(0,T;W^{2,\sigma}(\Omega)),\\
\label{regmu0}
&\mu\in L^\infty(0,T;V),\\
\label{regxi0}
&\xi\in  L^\infty(0,T;H)\cap L^\infty(0,T;L^\sigma(\Omega)),\\
\label{regw0}
&w\in H^2(0,T;H)\cap W^{1,\infty}(0,T;V),
\end{align}
where $\sigma$ is arbitrary in $(2,\infty)$, and we have the variational identities
\begin{align}
  & \iO \dt\phi \, v 
  + \iO \nabla\mu \cdot \nabla v
  + \gamma \iO \phi v
  = \iO f v
  \non
  \\	
  & \quad \hbox{for every $v\in V$ and a.e. $t\in (0,T)$}\,,
  \label{prima}
  \\
%    & \iO\mu v
%  = \iO\nabla\phi\cdot\nabla v + \iO\xi v + \iO (F'(\phi) + a - b \dt w)v
%	\nonumber\\
%	&  \quad \hbox{for every $v\in V$ and a.e. $t\in (0,T)$}\,,
%  \label{seconda}
%  \\
   & \iO\mu v
  = \iO\nabla\phi\cdot\nabla v + \iO\xi v + \iO (F'(\phi) + a - b \dt w)v \ 
  {}\hbox{ for every $v\in V$}	\nonumber\\
	&  \quad \hbox{and a.e. $t\in (0,T)$}\,,
%  \nonumber
%  \\
	 \hbox{ with } \, \xi \in \partial I_{[-1,1]}(\phi) \, \hbox{ a.e. in  $Q$}\,,
  \label{seconda}
  \\
    & \iO \dtt w \, v 
  + \iO \nabla( \kuno \dt w + \kdue w) \cdot \nabla v
  + \lambda \iO \dt\phi \, v
  = \iO u v
  \non
  \\
  & \quad \hbox{for every $v\in V$ and a.e. $t\in (0,T)$}\,,
  \label{terza}
  \\
  & \phi(0) = \phi_0, 
  \quad 
  w(0) = w_0,
  \quad 
  \dt w(0) = w_1 \,.
  \label{cauchy}
\end{align}
Moreover, the solution components $\,\phi\,$ and $\,w\,$ are uniquely determined, that is, whenever $(\phi_i,\mu_i,\xi_i,w_i)$, $i=1,2$, are two such
solutions, then $\phi_1=\phi_2$ and $w_1=w_2$.
\end{theorem}
\begin{remark}{\rm
By continuous embedding, we have $\phi,w,\dt w\in C^0([0,T];H)$ so that the evaluations $\phi(0),w(0),\dt w(0)$ and $\phi(T), w(T),\dt w(T)$ are meaningful. Moreover,
since the solution components $\phi$ and $w$ are uniquely determined, the cost functional (\ref{cost}) is well defined on $\Uad$.  
Besides, let us remark that the exponent $\sigma$ appearing above is more general than $\sigma=6$ which was stated in  
\cite[Thm.~2.2]{CGSS3}. This is due to the fact that here we only focus on the two-dimensional case, where the continuous embedding $V \subset 
L^\sigma(\Omega)$ is true for any $\sigma \in (2,\infty)$ {(the case $\sigma \in [1,2]$ is already ensured by the regularities $\phi \in L^\infty(0,T;W)$ and 
$\xi\in L^\infty(0,T;H)$}) instead of the three-dimensional embedding $V \subset L^6(\Omega)$.}
\end{remark}

For the deep quench approximations, we have the following stronger result.
\begin{theorem}
\label{SSal}
Suppose that {\rm (\ref{Uad}), (\ref{uminmax})} and {\bf (A1)}--{\bf (A4)} are satisfied. 
Then the system {\rm (\ref{Iprima}), (\ref{Isecondanew}), (\ref{Iterza})--(\ref{Icauchy})} has for every $u\in\Uad$ and every $\alpha>0$
a unique solution $(\phial,\mual,\wal)$ such that
\begin{align}
  & \phial \in W^{1,\infty}(0,T;H)\cap H^1(0,T;W)\cap L^\infty(0,T;W^{2,\sigma}(\Omega)),
    \label{regphial}
  \\
  & \mual \in L^\infty(0,T; V)\cap\Liq,
  \label{regmual}
  \\
	&	h_\alpha'(\phial)\in L^\infty(Q),
	\label{regxial}
	\\
  & \wal \in H^2(0,T; H) \cap C^1([0,T]; V),
  \label{regwal}
\end{align}
for arbitrary $\sigma\in (1,\infty)$. Moreover, there exists a constant $K_1(\sigma)>0$, which depends only on the structure of the system, 
$\Omega$, $T$, 
 the norms of the data, and the choice of $\sigma\in [2,\infty)$,
such that 
\begin{align}
  & \|\phial\|_{H^1(0,T; V)  \cap L^\infty(0,T;W^{2,\sigma}(\Omega))}
  \,+\,\|\mual\|_{L^\infty(0,T; V)}  \non\\
  &\,\,+\, \|h_\alpha'(\phial)\|_{L^\infty(0,T;L^{\sigma}(\Omega))}
    \,+\, \|\wal\|_{H^2(0,T;H) \cap C^1([0,T];V)}\,  \leq\, K_1(\sigma)\,,
  \label{boundal1}
	\end{align}
whenever $\alpha\in (0,1]$ and $u\in\Uad$. In addition, for every $\alpha>0$ there holds the {\it strict separation property}, i.e.,
there exist constants $-1<r_*(\alpha)<r^*(\alpha)<1$, which depend only on the structure of the system, $\Omega$, $T$, and the norms
of the data, such that for every $u\in\Uad$ it holds
\begin{equation}
\label{separation}
r_*(\alpha)\le \phial(x,t) \le r^*(\alpha) \quad\forall\,(x,t)\in \overline Q.
\end{equation} 
\end{theorem}
{\em Proof.} \,\,Existence, uniqueness and the regularity properties (\ref{regphial})--(\ref{regwal}) of the solution follow directly from 
\cite[Thms.~2.1~and~2.4]{CGSS4}. Moreover, \cite[Thm.~2.4]{CGSS4} also yields the existence of 
constants  $-1<r_*(\alpha)\le r^*(\alpha)<1$ such
that the inequality in (\ref{separation}) holds true at least for almost every $(x,t)\in Q$. But since $H^1(0,T;W)$ is continuously
embedded in $C^0(\overline Q)$, we have $\phial\in C^0(\overline Q)$, so that the pointwise condition (\ref{separation}) is in fact valid. 

It remains to show the existence of a constant $K_1(\sigma)$ satisfying (\ref{boundal1}). To this end, we recall the proof of 
\cite[Thm.~2.5]{CGSS3} (cf.~also 
\cite[Thm.~2.1 and Rem.~2.3]{CGSS4}). The strategy employed there was to approximate the system (\ref{Iprima}), (\ref{Isecondanew}), 
(\ref{Iterza})--(\ref{Icauchy}) by replacing the nonlinearities $\hal$ by their Moreau--Yosida approximations $h_{\alpha,\eps}$ 
at the level $\eps>0$. For sufficiently small $\eps>0$, existence, uniqueness, and uniform estimates could be shown for the
approximating system at the level $\eps>0$. 
Now observe that we have, for every $\alpha\in (0,1]$ and every $\eps>0$,
\begin{align*}
& 0\,\le\,h_{\alpha,\eps}(r)\,\le\,\hal(r)\,\le\,h(r) \quad\forall \,r\in\erre,\\
& |h_{\alpha,\eps}'(r)|\,\le\,|h_{\alpha}'(r)|\,\le\,|h'(r)|  \quad\forall \,r\in (-1,1).
\end{align*}
A closer inspection of the estimates performed in the cited proofs now reveals that the above uniform estimates, combined with the boundedness
of $\Uad$ and the fact that $\alpha\in (0,1]$,  have the consequence that all of the bounds derived in the cited proofs for the 
approximations at the level $\eps>0$ can in fact
be made uniformly with respect to the choice of $\alpha\in (0,1]$. Since all these estimates are also uniform with respect to 
sufficiently small $\eps>0$, they persist under the passage to the limit as $\eps\searrow0$, thanks to the semicontinuity of norms. 
This concludes the proof. \qed
\begin{remark}
\label{Variational}
{\rm 
The above well-posedness result in fact refers to the natural variational form \eqref{prima}
of the homogeneous Neumann problem for equation (\ref{Iprima}),
due to the low regularity of $\mual$ specified in~(\ref{regmual}).
However, thanks to (\ref{regphial}), {\bf (A3)}, and the elliptic regularity theory, it is clear
that $\mual\in L^\infty(0,T;W)$ as well, so that we actually can write (\ref{Iprima}) in its strong~form.
}
\end{remark}

Let us conclude this section by collecting some useful tools that will be employed later on.
We make frequent use of the Young, Poincar\'e--Wirtinger and compactness inequalities:
\begin{align}
  & ab \,\leq\, \delta a^2 + \frac 1{4\delta} \, b^2
  \quad \hbox{for every $a,b\in\erre$ and $\delta>0$},
  \label{young}
  \\[2mm]
    & \|v\|_V\,  \leq \,C_\Omega \, \bigl( \|\nabla v\| + |\vbar| \bigr)
  \quad \hbox{for every $v\in V$},
  \label{poincare}
  \\[2mm]
  & \|v\| \,  \leq\, \delta \, \|\nabla v\| + C_{\Omega,\delta} \, \|v\|_*
  \quad \hbox{for every $v\in V$ and $\delta>0$},
  \label{compact}
\end{align}
where $C_\Omega$ depends only on~$\Omega$, $C_{\Omega,\delta}$ depends on~$\delta$, in addition,
and $\|\,\cdot\,\|_*$ is the norm in $V^*$ to be introduced below (see (\ref{normaVp})).

Next, we recall an important  tool which is commonly used when working with problems connected to the Cahn--Hilliard equation.
Consider the weak formulation of the Poisson equation $-\Delta z=\psi$
with homogeneous Neumann boundary conditions. 
Namely, for a given $\psi\in V^*$ (and not necessarily in~$H$), we consider the problem:
\begin{align}
  \hbox{Find} \quad 
	z \in V
  \quad \hbox{such that} \quad
  \iO \nabla z \cdot \nabla v
  = \langle \psi , v \rangle
  \quad \hbox{for every $v\in V$}.
  \label{neumann}
\end{align}
Since $\Omega$ is connected and regular, it is well known
 that the above problem admits a family of solutions $z$ if and only if $\psi$ has zero mean value; among the solutions $z$ there is only one with zero mean value.
Hence, we can introduce the associated solution operator $\cal N$, which turns out to be an isomorphism between the following spaces, by
\begin{align}
  & \calN: \mbox{dom}(\calN):=\{\psi\in V^*:\ \psibar=0\} \to \{z\in V:\ \overline z=0\},
  \quad 
  {\cal N}: \psi \mapsto z,
  \label{defN}
\end{align}
where $z$  is the unique solution to \eqref{neumann} satisfying $\overline z=0$.
Moreover, it follows that the formula
\begin{align}
  \|\psi\|_*^2 := \|\nabla\calN(\psi-\psibar)\|^2 + |\psibar|^2
  \quad \hbox{for every $\psi\in V^*$}
  \label{normaVp}
\end{align}
defines a Hilbert norm in $\Vp$ that is equivalent to the standard dual norm of $\Vp$.
From the above properties, one can obtain the following identities:
\begin{align}
  & \iO \nabla\calN\psi \cdot \nabla v
  = \langle \psi , v \rangle
  \quad  \hbox{for every $\psi\in\mbox{dom}(\calN)$, $v\in V$},
  \label{dadefN}
  \\
  & \langle \psi , \calN\zeta \rangle
  = \langle \zeta , \calN\psi \rangle
  \quad \hbox{for every $\psi,\zeta\in\mbox{dom}(\calN)$},
  \label{simmN}
  \\
  & \langle \psi , \calN\psi \rangle 
  = \iO |\nabla\calN\psi|^2
  = \|\psi\|_*^2
  \quad \hbox{for every $\psi\in\mbox{dom}(\calN)$},
  \label{danormaVp}
\end{align}
as well as 
\begin{align}
  &\int_{t_1}^{t_2} \!\!\langle \dt v(s) , \calN v(s) \rangle \mbox{d}s
  = \int_{t_1}^{t_2} \!\!\langle v(s) , \calN(\dt v(s)) \rangle {\rm  d}s
  \nonumber
  \\
  & = \frac 12 \, \|v(t_2)\|_*^2
  - \frac 12 \, \|v(t_1)\|_*^2\,,
  \label{propN} 
\end{align}
which holds for all $t_1,t_2\in[0,T]$ with $t_1\le t_2$ and every $v\in H^1(0,T;\Vp)$ having zero mean value.

Finally, without further reference later on, we are going to employ the following convention: the capital-case symbol $C$ is used to denote every constant
that depends only on the structural data of the problem such as
$\Omega$, $T$, $a$, $b$, $\kuno$, $\kdue$, $\gamma$, $\lambda$, the shape of the 
nonlinearities, and the norms of the involved functions. Therefore, its meaning may vary from line to line and even within the same line.
In addition, when a positive constant $\delta$ enters the computation, then the related symbol $C_\delta$, in place of a general $C$, denotes constants
that depend on $\delta$, in addition.

%%%%%%%%%%%%%%%%%%%%%%%%%%%%%%%%%%%%%%%%%%%%%%%%%%%%%%%%%%%%%%%%%%%%%

\section{Deep quench approximation of states and optimal controls}

We begin our analysis by proving a result that provides a qualitative comparison between the deep quench approximations associated with
different values of $\alpha>0$.

\begin{theorem}
\label{SSal12}
Suppose that \eqref{Uad}, \eqref{uminmax} and {\bf (A1)}--{\bf (A4)} are fulfilled, and let, for given $u\in\Uad$ and 
$0<\alpha_1<\alpha_2\le 1$, the solutions to the deep quench system \eqref{Iprima}, \eqref{Isecondanew}, \eqref{Iterza}--\eqref{Icauchy}
given by {Theorem \ref{SSal}}
be denoted by $(\phi_{\alpha_i},\mu_{\alpha_i},w_{\alpha_i})$, $i=1,2$. Then there is a constant $K_2>0$, which depends only on the data
of the system, such that 
\begin{align}
\label{alpha12}
&\|\phi_{\alpha_1}-\phi_{\alpha_2}\|_{C^0([0,T];V^*)
\cap L^2(0,T;V)}\,+\,\|w_{\alpha_1}-w_{\alpha_2}\|_{H^1(0,T;H)\cap C^0([0,T];V)}   
\nonumber\\
&\le\,K_2\left(\alpha_2-\alpha_1\right)^{1/2}\,.
\end{align}
\end{theorem}

\noindent{\em Proof.} \,\,
We set, for convenience, 
\begin{align*}
&\phi:=\phi_{\alpha_1}-\phi_{\alpha_2},\quad  \mu:=\mu_{\alpha_1}-\mu_{\alpha_2},\quad w:=w_{\alpha_1}-w_{\alpha_2},\\
&\rho_i:=F'(\phi_{\alpha_i})\,\mbox{ for $i=1,2$}, \quad \rho:=\rho_1-\rho_2\,.
\end{align*}
Then $(\phi,\mu,w)$ is a solution to the system which in its strong formulation reads as follows:
\begin{alignat}{2}
\label{diff1}
& \dt\phi-\Delta\mu + \gamma \phi=0 
&&\quad\mbox{in $Q$},\\
\label{diff2}
& \mu =-\Delta\phi +h_{\alpha_1}'(\phi_{\alpha_1})-h_{\alpha_2}'(\phi_{\alpha_2})+\rho-b\dt w
&&\quad\mbox{in $Q$},\\
\label{diff3}
&\dtt w -\Delta(\kappa_1\dt w+\kappa_2 w)+\lambda \dt\phi=0
&&\quad\mbox{in $Q$},\\
\label{diff4}
&\dn\phi(0)=\dn \mu=\dn(\kappa_1 \dt w+ \kappa_2 w)=0 
&&\quad\mbox{on $\Sigma$},\\
\label{diff5}
&\phi(0)=w(0)=\dt w(0)=0 &&\quad\mbox{in  $\Omega$}.
\end{alignat}

We first observe that the mean values of $\phi$ and $\dt\phi$ vanish on $[0,T]$. Indeed, testing \eqref{diff1} by the constant
function $\boldsymbol{1}/|\Omega|$ yields
that
\begin{equation}
\frac d{dt} \,\overline\phi(t)+\gamma\overline\phi(t)=0 \quad\forall \,t\in(0,T], \quad \overline\phi(0)=0\,,
\end{equation}
whence the claim readily follows. Therefore, the expression ${\cal N}\phi$ is meaningful as an element of $V$. 
We now test \eqref{diff1} by ${\cal N}\phi$, \eqref{diff2} by $\phi$, and we integrate \eqref{diff3} with respect time over
$[0,t]$ and test the resulting identity by $\,\frac b \lambda\,\dt w$. Then we add the three resulting equations 
to each other. Using the properties \eqref{dadefN}--\eqref{propN}, we find that four terms cancel, and it follows
the identity
\begin{align}
\label{Lollouno}
&\frac 12 \,\|\phi(t)\|_*^2\,+\,\gamma\int_0^t \|\phi(s)\|_*^2 {\rm d}s +\intQt |\nabla\phi|^2\,+\,\frac b \lambda \intQt |\dt w|^2
\non\\
&\quad +\,\frac {b\kappa_1}{2\lambda}\iO|\nabla w(t)|^2 \,+\intQt 
\bigl(h_{\alpha_1}'(\phi_{\alpha_1})-h_{\alpha_1}'(\phi_{\alpha_2})\bigr)\,\phi \nonumber\\
&=\,-\intQt \bigl(h_{\alpha_1}'(\phi_{\alpha_2})-h_{\alpha_2}'(\phi_{\alpha_2})\bigr)\,\phi\,-\intQt \rho\,\phi
\,-\,\frac{b \kappa_2}\lambda\intQt  (\boldsymbol{1} * \nabla w) \cdot \nabla\dt w  \nonumber\\
&=:\, I_1+I_2+I_3\,,
\end{align}
with obvious meaning. Owing to the monotonicity of $\,h_{\alpha_1}'$, the last term on the left-hand side is nonnegative. Moreover, thanks to
the fact that $\phi_{\alpha_1}$ and $\phi_{\alpha_2}$ attain their values in $(-1,1)$, it follows from the convexity of $\,h\,$ that
\begin{align*}
&-\,\bigl(h_{\alpha_1}'(\phi_{\alpha_2})-h_{\alpha_2}'(\phi_{\alpha_2})\bigr)\,\phi
\,=\,(\alpha_2-\alpha_1)\,h'(\phi_{\alpha_2})\,\phi \nonumber\\
&\le\,(\alpha_2-\alpha_1)(h(\phi_{\alpha_1})-h(\phi_{\alpha_2}))\,
\le\,(\alpha_2-\alpha_1)\,2\,\ln(2)\,,
\end{align*}
so that 
\begin{equation}
\label{Lollodue}
I_1\,\le\,(\alpha_2-\alpha_1)\,2\,\ln(2)\,|\Omega|\,T\,.
\end{equation}
Moreover, invoking the Lipschitz continuity of $F'$, as well as the compactness inequality \eqref{compact}, we conclude that
\begin{align}
\label{Lollotre}
|I_2| &\,\le\,\intQt \left|F'(\phi_{\alpha_1})-F'(\phi_{\alpha_2})\right|\,|\phi|\,\le\,C\intQt |\phi|^2\nonumber\\
&\,\le\frac 12\intQt |\nabla\phi|^2 \,+\,C\int_0^t \|\phi(s)\|_*^2 {\rm d}s\,.
\end{align}
It remains to estimate $I_3$. To this end, using the identity
\begin{align*}
	\intQt (\boldsymbol{1}*\nabla w) \cdot \nabla (\dt w )
	= \iO (\boldsymbol{1}*\nabla w(t)) \cdot \nabla w(t)
	- \intQt  |\nabla w|^2, 
\end{align*}
the fact that 
\,$\Vert \boldsymbol{1}*\nabla w(t) \Vert^2 \leq \Big( \int_0^t \Vert \nabla w (s)
\Vert {\rm d}s\Big)^2\leq T\intQt |\nabla w|^2$, 
as well as Young's inequality, we infer that
\begin{align}
\label{Lolloquattro}
	I_3 
	\,\leq\, \frac{b\kappa_1}{4\,\lambda}\, \|\nabla w(t)\|^2
	+ C \intQt |\nabla w|^2.
\end{align}
Combining \eqref{Lollouno}--\eqref{Lolloquattro}, and invoking Gronwall's lemma, 
we have thus shown that 
\begin{align*}
&\|\phi\|_{L^\infty(0,T;V^*)}^2\,
+\, \|\nabla \phi\|^2_{L^2(0,T; H)}
\,
+\,\|w\|_{H^1(0,T;H)\cap L^\infty(0,T;V)}^2
%\nonumber\\
%&
\le\,C\,(\alpha_2-\alpha_1)\,.
\end{align*}
The assertion now follows from the fact that the $L^2(Q)$ norm of $\,\phi\,$ can be estimated via the compactness inequality \eqref{compact}. 
\qed
 
\begin{theorem}
\label{THM:anktozero}
Suppose that \eqref{Uad}, \eqref{uminmax} and {\bf (A1)}--{\bf (A4)} are fulfilled, and let sequences $\{\alpha_n\}\subset (0,1]$
and $\{u_{\alpha_n}\}\subset\Uad$ be given such that $\alpha_n\searrow0$ and $u_{\alpha_n}\to u$ weakly star in $\Liq$ for some $u\in\Uad$. 
Moreover, let 
$(\phi_{\alpha_n},\mu_{\alpha_n},w_{\alpha_n})$ be the solution in the sense of {Theorem \ref{SSal}} to the problem 
{\eqref{Iprima}, \eqref{Isecondanew}, \eqref{Iterza}--\eqref{Icauchy}} with the control $u_{\alpha_n}$ and the convex function $h_{\alpha_n}$,
for $n\in\enne$. Then there exist a subsequence $\{\alpha_{n_k}\}$ and a solution $(\phi^0,\mu^0,\xi^0,w^0)$ 
to the state system \eqref{Iprima}--\eqref{Icauchy}
such that, as $k\to\infty$,
\begin{alignat}{2}
\label{conphial}
\phi_{\alpha_{n_k}} & \to \phi^0 && \quad\mbox{weakly star in }\, H^1(0,T;V)\cap L^\infty(0,T; W^{2,\sigma}(\Omega))
\nonumber
\\
& && \quad \mbox{and strongly in }\,C^0(\overline Q),\\
\label{conmual}
\mu_{\alpha_{n_k}}& \to \mu^0 && \quad\mbox{weakly star in }\,L^\infty(0,T;V),\\
\label{conxial}
h_{\alpha_{n_k}}'(\phi_{\alpha_{n_k}})& \to \xi^0 &&\quad\mbox{weakly star in }\,L^\infty(0,T; L^{\sigma}(\Omega)), \\
\label{conwal}
w_{\alpha_{n_k}} & \to w^0 && \quad\mbox{weakly star in }\,H^2(0,T;H)\cap W^{1,\infty}(0,T;V)
\nonumber
\\
& &&\quad \mbox{and strongly in }\,C^1([0,T];H),
\end{alignat}
with $\sigma$ arbitrary in $[2,\infty)$.
\end{theorem}
\noindent{\em Proof.} \,\,By virtue of the global estimate \eqref{boundal1}, it follows the existence of 
the subsequence and of limits $(\phi^0,\mu^0,\xi^0,w^0)$
such that \eqref{conphial}--\eqref{conwal} hold true. In this connection, the strong convergence result in \eqref{conphial} 
follows from standard compactness results (see, e.g., \cite[Sect.~8, Cor.~4]{Simon}). Observe that
the strong convergence in \eqref{conphial} and the Lipschitz continuity of $F'$ imply that 
$F'(\phi_{\alpha_{n_k}}) \to F'(\phi^0)$ strongly in $C^0(\overline Q)$ {as $k \to \infty$}. 

We then need to show that $(\phi^0,\mu^0,\xi^0,w^0)$ is a solution to \eqref{Iprima}--\eqref{Icauchy}. Owing to the convergence properties
\eqref{conphial}--\eqref{conwal}, it is easily verified by passage to the limit as $k\to\infty$ that $(\phi^0,\mu^0,\xi^0,w^0)$ satisfies the 
(equivalent) time-integrated version of the variational equalities in~\eqref{prima}--\eqref{terza} with test functions $v\in L^2(0,T;V)$ for the control $u$.
Also, the initial conditions in \eqref{cauchy} follow easily from the weak convergences in \eqref{conphial} and \eqref{conwal}.
It remains to show that $\,\xi^0\in\sindi(\phi^0)$ almost everywhere in $Q$. For this purpose, we define on $L^2(Q)$ the convex functional
\begin{equation*}
\Phi(v)=\int_Q \indi(v), \,\mbox{ if }\indi(v)\in L^1(Q),\,\mbox{ and }\, \Phi(v)=+\infty,\,\mbox{ otherwise.}
\end{equation*}
It then suffices to show that $\xi^0$ belongs to the subdifferential of $\Phi$ at $\phi^0$, i.e., that
\begin{equation}
\label{Egon}
\Phi(v)-\Phi(\phi^0)\,\ge\,\int_Q \xi^0(v-\phi^0)\quad\forall\,v\in L^2(Q).
\end{equation}
At this point, recall that $\phi_{\alpha_{n_k}}(x,t)\in [-1,1]$, and thus also {$\phi^0(x,t)\in [-1,1]$ in $\overline Q$.}
Consequently, $\Phi(\phi^0)=0$. Now observe that in the case that $\indi(v)\not\in L^1(Q)$ the inequality \eqref{Egon} holds true since 
its left-hand side is infinite. If, however, $\indi(v)\in L^1(Q)$, then obviously $v\in [-1,1]$ almost everywhere in $Q$, and it follows from
\eqref{halcon} and Lebesgue's theorem of dominated convergence that
$$
\lim_{k\to\infty}\,\int_Q h_{\alpha_{n_k}}(v)\,=\,\Phi(v)\,=\,0.
$$
Now, by the convexity of $h_{\alpha_{n_k}}$, and since $h_{\alpha_{n_k}}(\phi_{\alpha_{n_k}})$ is nonnegative, we have for all $v\in L^2(Q)$
that
\begin{align*}
h_{\alpha_{n_k}}'(\phi_{\alpha_{n_k}})(v-\phi_{\alpha_{n_k}})\,\le\,h_{\alpha_{n_k}}(v)-h_{\alpha_{n_k}}(\phi_{\alpha_{n_k}})\,\le\,h_{\alpha_{n_k}}(v)
\quad\mbox{a.e. in }\,Q.
\end{align*}
Using \eqref{conphial} and \eqref{conxial}, we thus obtain the following chain of (in)equalities: 
\begin{align*}
\int_Q \xi^0(v-\phi^0) \,&=\,\lim_{k\to\infty} \int_Q h'_{\alpha_{n_k}}(\phi_{\alpha_{n_k}})(v-\phi_{\alpha_{n_k}}) 
\\
&\le\,\limsup_{k\to\infty}\int_Q \big(h_{\alpha_{n_k}}(v)-h_{\alpha_{n_k}}(\phi_{\alpha_{n_k}})\big)
\\
&\leq\, \lim_{k\to\infty} \int_Q h_{\alpha_{n_k}}(v) \,=\,\Phi(v)\,=\,\Phi(v)-\Phi(\phi^0),
\end{align*}
which shows the validity of \eqref{Egon}. This concludes the proof.
\qed
\begin{remark} \label{RMK:limit}
{\rm Since, according to {Theorem \ref{SS0}}, the solution variables $\phi^0$ and $w^0$ are uniquely determined, the convergence properties
\eqref{conphial} and \eqref{conwal} actually hold for the entire sequences and not just for a subsequence.} 
%Moreover, the above theorem shows that the state system \eqref{Iprima}--\eqref{Icauchy} 
%admits at least one solution $(\phi,\mu,\xi,w)=(\phi^0,\mu^0,\xi^0,w^0)$ that enjoys
%the same regularity properties as the deep quench approximations, namely 
%\eqref{regphial}--\eqref{regwal}. A fortiori, since the 
%solution  components $\phi^0$ and $w^0$ are uniquely determined, the corresponding 
%solution components have this regularity for any solution.
\end{remark}
\begin{corollary}
Suppose that  \eqref{Uad}, \eqref{uminmax} and {\bf (A1)}--{\bf (A4)} are fulfilled, let $(\phi^0,\mu^0, \xi^0,w^0)$
be a solution to \eqref{Iprima}--\eqref{Icauchy} and $(\phial,\mual,\wal)$ be the solution to 
{\eqref{Iprima}, \eqref{Isecondanew}, \eqref{Iterza}--\eqref{Icauchy}} associated with some $\alpha\in (0,1]$. 
Then, with the constant $K_2>0$ introduced in Theorem {\ref{SSal12}}, we have  
\begin{align}
\label{messi}
&\|\phi_\alpha-\phi^0\|_{C^0([0,T];V^*)
\cap L^2(0,T;V)}\,+\,\|w_{\alpha}-w^0\|_{H^1(0,T;H)\cap C^0([0,T];V)}   
%\nonumber\\ &
\le\,K_2\,\alpha^{1/2}\,.
\end{align}
\end{corollary}
\noindent{\em Proof.}
This is an immediate consequence of the uniqueness of $w^0$ and $\phi^0$, if we put in \eqref{alpha12} $\,\alpha_2=\alpha\,$ and pass to
the limit as $\alpha_1\searrow0$.
\qed
	
%%%%%%%%%%%%%%%%%%%%%%%%%%%%%%%%%%%%%%%%%%%%%%%%%%%%%%%%%%%%%%%%%

\section{Existence and approximation of optimal controls}

Beginning with this section, we study the optimal control problem  {\bf (P)} of minimizing the cost functional \eqref{cost} subject to the state 
system \eqref{Iprima}--\eqref{Icauchy} and the control constraint $u\in\Uad$, where \eqref{Uad} and \eqref{uminmax} are generally assumed
to be valid. In addition to
{\bf (A1)}--{\bf (A4)}, we impose the following general assumptions:

\begin{description}
\item[(A5)] \,\,The coeffients $\beta_1,\ldots,\beta_6,\nu$ are nonnegative and not all equal to zero.
\item[(A6)] \,\,$\phi_\Omega,w_\Omega,w_\Omega' \in \Ldue$\, and \,$\phi_Q,w_Q,w_Q'\in L^2(Q)$.      
\end{description}
We compare the problem {\bf (P)} with the following family of optimal control problems for $\alpha>0$:

\vspace{2mm}\noindent
{\bf (P$_\alpha$)}  \,\,\,Minimize the cost functional \eqref{cost} subject to the state system 
\eqref{Iprima}, \eqref{Isecondanew}, \hspace*{11.5mm}\eqref{Iterza}--\eqref{Icauchy} and the control constraint $u\in\Uad$.

\vspace{2mm}\noindent
We expect that the minimizers of the control problems {\bf (P)} and {\bf (P$_\alpha$)} are closely related. Before giving an affirmative answer
to this conjecture, we introduce for convenience the following control-to-state operators:
\begin{align}
{\cal S}:\, &\Uad\ni u\mapsto (\phi,w), \,\mbox{ where $\phi,w$ are the first and fourth components}
\nonumber\\
&\mbox{of any solution to 
\eqref{Iprima}--\eqref{Icauchy}}, 
\label{defS}
\\
{\cal S}_\alpha :\, &\Uad\ni u\mapsto (\phial,\wal), \,\mbox{ where $\phial,\wal$ are the first and third components}
\nonumber\\
&\mbox{of the solution 
to {\eqref{Iprima}, \eqref{Isecondanew}, \eqref{Iterza}--\eqref{Icauchy}}}.
\label{defSal}
\end{align}
We then have the following result.          
\begin{proposition}
\label{PROP:J}
Suppose that \eqref{Uad}, \eqref{uminmax} and {\bf (A1)}--{\bf (A6)} are fulfilled, and let sequences $\{\alpha_n\}\subset
(0,1]$ and $\{u_n\}\subset\Uad$ be given such that $\alpha_n\searrow0$ and $u_n\to u$ weakly star in $\Liq$ for some $u\in\Uad$ as $n \to\infty$.
Then we have
\begin{align}
\label{uwe1}
\J(\S(u),u)\,&\le\,\liminf_{n\to\infty}\, \J({\cal S}_{\alpha_n}(u_n),u_n),\\
\label{uwe2}
\J(\S(v),v)\,&=\,\lim_{n\to\infty} \, \J({\cal S}_{\alpha_n}(v),v) \quad\forall\,v\in\Uad.
\end{align}
\end{proposition} 
{\em Proof.} \,\,According to Theorem {\ref{THM:anktozero}} and Remark {\ref{RMK:limit}}, 
we have $\Saln(u_{\alpha_n})\to \S(u)$ in the sense of \eqref{conphial} and \eqref{conwal},
respectively. Then \eqref{uwe1} follows from the semicontinuity properties of the cost functional. Now let $v\in\Uad$ be 
arbitrarily chosen. Applying Theorem {\ref{THM:anktozero}} and Remark \ref{RMK:limit} to the constant sequence $v_n=v$, $n\in\enne$, we infer
that $\Saln(v)\to \S(v)$ in the sense of \eqref{conphial} and \eqref{conwal}. In particular, this implies strong convergence of the
sequences $\{\phialn\}$, $\{\waln\}$ and $\{\dt \waln\}$ in $C^0([0,T];H)$, by compact embedding. 
Since the first six summands of the cost functional are
continuous with respect to the strong topology of $C^0([0,T];H)$, the validity of \eqref{uwe2} follows.
\qed

\vspace{2mm}
We are now in a position to show the existence of minimizers for the control problem {\bf (P)}. We have the following result.
\begin{corollary}
\label{COR:exoptimal}
Suppose that \eqref{Uad}, \eqref{uminmax}, and {\bf (A1)}--{\bf (A6)} are fulfilled. Then the problem {\bf (P)}
admits at least one solution in $\Uad$.
\end{corollary}
{\em Proof.} \,\,We pick an arbitrary sequence $\{\alpha_n\}\subset (0,1]$ such that $\alpha_n\searrow0$ as $n\to\infty$. By virtue of 
\cite[Thm.~4.1]{CGSS3}, the problem {\bf (P$_{\alpha_n}$)} has a solution $\,u_{\alpha_n}\in\Uad\,$ with associated state 
$(\phialn,\mualn,\waln)$ and $\xi_{\alpha_n}:=\xialn$ for  $n\in\enne$. Since $\Uad$ is bounded in $\Liq$, we may without loss
of generality assume that $u_{\alpha_n}\to u$ weakly star in $\Liq$ for some $u\in\Uad$. Then, in view of Theorem \ref{THM:anktozero},
there are a subsequence $\{\alpha_{n_k}\}$ and a solution $(\phi^0,\mu^0,\xi^0,w^0)$ to the system \eqref{Iprima}--\eqref{Icauchy} such that the 
convergence properties \eqref{conphial}--\eqref{conwal} hold true. Now observe that $(\phi^0,w^0)=\S(u)$ and $(\phialn,\waln)=\Saln(u_n)$
for $n\in\enne$. We then obtain from the optimality of $((\phialn,\waln),u_{\alpha_n})$ for {\bf (P$_{\alpha_n}$)}, using Proposition \ref{PROP:J},
the following chain of (in)equalities:
\begin{align*}
\J(\S(u),u)\,&\le\,\liminf_{k\to\infty} \,\J({\cal S}_{\alpha_{n_k}}(u_{\alpha_{n_k}}),u_{\alpha_{n_k}})\,\le
\, 	\liminf_{k\to\infty}\,\J({\cal S}_{\alpha_{n_k}}(v),v) \nonumber\\
&=\,\J(\S (v),v)\,.
\end{align*}
This shows that $(\S(u),u)$ is an optimal pair of the control problem {\bf (P)}, which concludes the proof 
of the assertion.
\qed

\vspace{2mm}
The proof of Corollary~\ref{COR:exoptimal} suggests that optimal controls of {\bf (P$_\alpha$)}
are ``close'' to optimal
controls of {\bf (P)} as $\alpha$ approaches zero. However, they do not yield any information on 
whether every optimal control of {\bf (P)} can be approximated in this way. In fact, such a global result cannot be expected to hold true. 
Nevertheless, a local answer can be given by employing a well-known trick. To this end, 
let $u^*\in\Uad$ be an optimal control for {\bf (P)} and $(\phi^*,\mu^*,\xi^*,w^*)$ be a solution to the
associated state system \eqref{Iprima}--\eqref{Icauchy} so that $(\phi^*,w^*)=\S(u^*)$. We associate
with this optimal control the {\em adapted cost functional}
\begin{equation}
\label{adcost}
\widetilde{\cal J}((\phi,w),u):=
{\cal J}((\phi,w),u)\,+\,\frac 12\,\|u-u^*\|^2_{L^2(Q)}
\end{equation}
and a corresponding \emph{adapted optimal control problem} for $\alpha>0$, namely:

\vspace{3mm}\noindent
 {\bf (${\widetilde{\bf P}}_{\alpha}$)}
 \quad Minimize $\,\, 
\widetilde {\cal J}((\vp,w),u)\,\,$
for $\,u\in\Uad$\,  subject to \,$(\vp,w)=\Sal
(u)$.

\vspace{3mm}\noindent
With essentially the same proof as that of \cite[Thm.~4.1]{CGSS4} (which needs no repetition here), we can show that the 
adapted optimal control problem 
{ \bf (${\widetilde{\bf P}}_{\alpha}$)} has for every $\alpha>0$ at least one solution.
The following result gives a partial answer to the question raised above
concerning the approximation of optimal controls for {\bf (P)} by the approximating problem {\bf (${\widetilde{\bf P}}_{\alpha}$)}.

\begin{theorem} \label{THM:APPROX}
Suppose that \eqref{Uad}, \eqref{uminmax} and {\bf (A1)}--{\bf (A6)} are fulfilled, assume that
$u^*\in \Uad$ is an arbitrary optimal control of {\bf (P)} with associated state  
$(\phi^*,\mu^*,\xi^*, w^*)$, and let $\,\{\alpha_k\}_{k\in\enne}\subset 
(0,1]\,$ be
any sequence such that $\,\alpha_k\searrow 0\,$ as $\,k\to\infty$. 
Then, for any $k\in \enne$, there exists 
an optimal control
 $\,u_{\alpha_{k}}\in \Uad\,$ of the adapted problem 
{ \bf (${\widetilde{\bf P}}_{\alpha_k}$)}
 with associated state $(\vp_{\alpha_k},\mu_{\alpha_k},w_{\alpha_k})$
 such that, as $k\to\infty$,
\begin{align}
\label{tr3.4}
&u_{\alpha_{k}}\to u^* \quad\mbox{strongly in }\,L^2(Q),
\end{align}
and such that \eqref{conphial}--\eqref{conwal} hold true with some $(\vp^0,\mu^0,\xi^0,w^0)$ satisfying 
$\phi^0=\phi^*$ and $w^0=w^*$. 
 Moreover, we have 
\begin{align}
\label{tr3.5}
&\lim_{k\to\infty}\,\widetilde{{\cal J}}(\S_{\alpha_{k}}(u_{\alpha_{k}}),u_{\alpha_{k}})
\,=\,{\cal J}(\S(u^*),u^*).
\end{align}
\end{theorem}
{\em Proof.} \,\,For any $ k\in\enne$, we pick an optimal control
$u_{\alpha_k} \in \Uad\,$ for the adapted problem { \bf (${\widetilde{\bf P}}_{\alpha_k}$)} and denote by 
$(\vp_{\alpha_k},\mu_{\alpha_k},w_{\alpha_k})$ the associated strong solution to the approximating
state system {\eqref{Iprima}, \eqref{Isecondanew}, \eqref{Iterza}--\eqref{Icauchy}}.
By the boundedness of $\Uad$ in $\Liq$, there is some subsequence $\{\alpha_{n}\}$ 
of $\{\alpha_k\}$ such that
\begin{equation}
\label{ugam}
u_{\alpha_{n}}\to u\quad\mbox{weakly star in }\,\Liq
\quad\mbox{as }\,n\to\infty,
\end{equation}
for some $u\in\Uad$. 
Thanks to Theorem~{\ref{THM:anktozero}}, the convergence properties \eqref{conphial}--\eqref{conwal} hold true correspondingly
for some solution $(\vp^0,\mu^0,\xi^0,w^0)$ to the state system \eqref{Iprima}--\eqref{Icauchy}, 
and  the pair $(\S(u),u)=((\vp^0,w^0),u)$
is admissible for {\bf (P)}. 

We now aim at showing that $ u=u^*$. Once this is shown, it follows from the uniqueness of the first and fourth components 
of the solutions to the state system \eqref{Iprima}--\eqref{Icauchy}
that also $(\vp^0,w^0)=
(\phi^*,w^*)$. 
Now observe that, owing to the weak sequential lower semicontinuity properties of 
$\widetilde {\cal J}$, 
and in view of the optimality property of $(\S(u^*),u^*)$ 
for problem {\bf (P)},
\begin{align}
\label{tr3.6}
\liminf_{n\to\infty}\, \widetilde{\cal J}(\S_{\alpha_n}(u_{\alpha_n}),u_{\alpha_n})\,
&\ge \,{\cal J}(\S(u),u)\,+\,\frac{1}{2}\,
\|u-u^*\|^2_{L^2(Q)}\nonumber\\[1mm]
&\geq \, {\cal J}(\S(u^*),u^*)\,+\,\frac{1}{2}\,\|u-u^*\|^2_{L^2(Q)}\,.
\end{align}
On the other hand, the optimality property of  
$\,(\S_{\alpha_{n}}(u_{\alpha_n}),u_{\alpha_n})
\,$ for problem { \bf (${\widetilde{\bf P}}_{\alpha_n}$)} 
yields that
for any $n\in\enne$ we have
\begin{equation}
\label{tr3.7}
\widetilde {\cal J}({\cal S}_{\alpha_{n}}(u_{\alpha_{n}}),
u_{\alpha_{n}})\,\le\,\widetilde {\cal J}({\cal S}_{\alpha_n}
(u^*),u^*)\,=\, {\cal J}({\cal S}_{\alpha_n}
(u^*),u^*)\,,
\end{equation}
whence, taking the limit superior as $n\to\infty$ on both sides and invoking (\ref{uwe2}) in
Proposition~{\ref{PROP:J}},
\begin{align}
	\non
	&\limsup_{n\to\infty}\,\widetilde {\cal J}(\S_{\alpha_{n}}(u_{\alpha_n}),
u_{\alpha_n})\,\le\,\limsup_{n\to\infty}\widetilde {\cal J}(\S_{\alpha_n}(u^*),u^*) 
\\ &  
\label{tr3.8}
=\,\limsup_{n\to\infty} {\cal J}(\S_{\alpha_n}(u^*),u^*)
\,=\,{\cal J}(\S(u^*),u^*)\,.
\end{align}
Combining (\ref{tr3.6}) with (\ref{tr3.8}), we have thus shown that 
$\,\frac{1}{2}\,\|u-u^*\|^2_{L^2(Q)}=0$\,,
so that $\,u=u^*\,$  and thus also $(\phi^*,w^*)
=(\vp^0,w^0)$. 
Moreover, (\ref{tr3.6}) and (\ref{tr3.8}) also imply that
\begin{align*}
%\label{tr3.9}
&{\cal J}(\S(u^*),u^*) \, =\,\widetilde{\cal J}(\S(u^*),u^*)
\,=\,\liminf_{n\to\infty}\, \widetilde{\cal J}(\S_{\alpha_n}(u_{\alpha_n}),
 u_{\alpha_{n}})\nonumber\\[1mm]
&=\,\limsup_{n\to\infty}\,\widetilde{\cal J}(\S_{\alpha_n}(u_{\alpha_n}),
{u_{\alpha_n}}) \,
=\,\lim_{n\to\infty}\, \widetilde{\cal J}(\S_{\alpha_n}(u_{\alpha_n}),
u_{\alpha_n})\,,
\end{align*}                                    
which proves the validity of {(\ref{tr3.5})}. 
Moreover, the convergence properties 
\eqref{conphial}--\eqref{conwal} are satisfied with $\phi^0=\phi^*$ and $w^0=w^*$. On the other hand, we have that 
\begin{align*}
%\label{tr3.9bis}
{\cal J}(\S(u^*),u^*) 
\,&\leq\,\liminf_{n\to\infty}\, {\cal J}(\S_{\alpha_n}(u_{\alpha_n}),
 u_{\alpha_{n}})
 \,\leq\,\limsup_{n\to\infty}\, {\cal J}(\S_{\alpha_n}(u_{\alpha_n}),
 u_{\alpha_{n}}) \nonumber\\[1mm]
&\leq \,\limsup_{n\to\infty}\,\widetilde{\cal J}(\S_{\alpha_n}(u_{\alpha_n}),
{u_{\alpha_n}}) \,
=\,{\cal J}(\S(u^*),u^*) ,
\end{align*}        
so that also 
$ {\cal J}(\S_{\alpha_n}(u_{\alpha_n}),{u_{\alpha_n}})$ converges to $ {\cal J}(\S(u^*),u^*)$ as 
$n\to \infty $, and the relation in \eqref{adcost} enables us to infer the strong convergence in \eqref{tr3.4}
for the subsequence $\{u_{\alpha_n}\}$. 

We now claim that \eqref{tr3.4} and \eqref{tr3.5} hold true even for the entire sequence, due to the complete identification of the limit
$u$ as $u^*$. We only prove this claim for \eqref{tr3.4}; for \eqref{tr3.5} a similar reasoning may be used.
Assume that \eqref{tr3.4} were not true. Then there exist some $\varepsilon>0$ and a subsequence $\{\alpha_j\}$ of  $\{\alpha_k\}$ such that  
\begin{equation}
\label{sirtoby}
\|u_{\alpha_j}-u^*\|_{L^2(Q)}\,\ge\,\varepsilon \quad\forall\,j\in\enne.
\end{equation}
However, by the boundedness of $\Uad$, there is some subsequence $\{\alpha_{j_n}\}$ 
of $\{\alpha_j\}$ such that, with some $\tilde u\in\Uad$,
\begin{equation*}
u_{\alpha_{j_n}}\to \tilde u\quad\mbox{weakly star in }\,\Liq
\quad\mbox{as }\,n\to\infty\,.
\end{equation*}
Arguing as above, it then turns out that $\tilde u = u^* $ and that \eqref{tr3.4} holds for the subsequence
$\{u_{\alpha_{j_n}}\}$ as well, which contradicts the fact that \eqref{sirtoby} obviously implies
that $\{ u_{\alpha_{j}} \}$ cannot have a subsequence that
converges strongly to $\,u^*\,$ in $L^2(Q)$.
\qed

%%%%%%%%%%%%%%%%%%%%%%%%%%%%%%%%%%%%%%%%%%%%%%%%%%%%%%%%%%%%%%%%%%%%%%%%%%%

\section{First-order necessary optimality conditions}

We now derive first-order necessary optimality conditions for
the control problem {\bf (P)}, using the corresponding conditions
for {\bf (${\widetilde{\bf P}}_{\alpha}$)} as approximations. To this end, we generally assume that the conditions \eqref{Uad}, 
\eqref{uminmax}, and {\bf (A1)}--{\bf (A6)} are fulfilled. Moreover, we need an additional assumption:

\begin{description} 
\item[(A7)]  \,\,It holds that $\,\beta_2\, \varphi_\Omega\in V\,$ and $\,\beta_6 \, w'_\Omega\in V$.
\end{description}
Notice that this assumption essentially requires a better regularity for the target data $
\varphi_\Omega$ and $\,w'_\Omega$ that in the cost functional give the endpoint tracking 
for the variables $\phi$ and $\dt w$. On the other hand, the assumption~{\bf (A7)} 
is trivially satisfied if $\beta_2 = \beta_6=0$.

Now let $u^*\in\Uad$ be a fixed optimal control for {\bf (P)} with
associated state $(\phi^*,\mu^*,\xi^*, w^*)$ (where only $\phi^*$ and $w^*$ are uniquely 
determined), and assume that $\alpha\in(0,1]$ is fixed. 
Moreover, suppose that $u_\alpha^*\in\Uad$ is an optimal control for {\bf (${\widetilde{\bf P}}_{\alpha}$)} with corresponding state
$(\phi_\alpha^*,\mu_\alpha^*,w_\alpha^*)$. The corresponding adjoint problem is given, in its strong form for simplicity, by
\begin{align}
	& - \dt \pal - \Delta \qal +\gamma\pal + \hal''(\phi_\alpha^*) \qal+F''(\phi_\alpha^*)\qal -\lambda\dt\ral
	\nonumber\\
	&\quad =\,\beta_1(\phi_\alpha^*-\phi_Q) \quad \mbox{in }\,Q,
	\label{adsys1}\\
	& 	 \qal \,=\, -  \Delta \pal \quad\mbox{in }\,Q,
	\label{adsys2}
	\\ 
	& \nonumber
	-\dt\ral -\Delta(\kappa_1\ral-\kappa_2(\boldsymbol{1}\circledast \ral)) -b\qal \\
	&\quad =\, \beta_3 ( \boldsymbol{1} \circledast (w_\alpha^* - w_Q) )+\beta_4(w^*_\alpha(T)-w_\Omega)
	+ \beta_5 (\dt w^*_\alpha - w_Q') \quad\mbox{in }\,Q,
	\label{adsys3}
	\\
& \dn \pal= \dn \qal = \dn (\kappa_1\ral -\kappa_2(\boldsymbol{1}\circledast \ral))=0  \quad\mbox{on }\,\Sigma,
\label{adsys4}
\\%[2mm]
& \pal(T)=  \beta_2 (\phi_\alpha^*(T)- \phi_\Omega) - \lambda{\beta_6} (\dt \wal^* (T)- w_\Omega') 
\nonumber\\
&\quad{}\hbox{and}\quad  \ral(T)=\beta_6 (\dt \wal^*(T)- w_\Omega')
\quad\mbox{in }\,\Omega,
\label{adsys5}
\end{align}
with the convolution product $\circledast$ defined in \eqref{conv:back}. 
Concerning this product, please note that $\dt (1 \circledast r) = -r.$
Let us, for convenience, introduce the abbreviations
\begin{align}
\label{source}
	&f_{\alpha}:= \beta_3 ( \boldsymbol{1} \circledast (w^*_\alpha - w_Q) )
	+ \beta_5 (\dt w^*_\alpha - w_Q')
	+ \beta_4 (w^*_\alpha(T) - w_\Omega),\\
\label{g1al}
&g_{\alpha}:=	\beta_1(\phi_\alpha^*-\phi_Q) ,\\
\label{rhoal}
&\rho_{\alpha}:=\beta_6 (\dt \wal^*(T)- w_\Omega') ,\\
\label{pial}
&\pi_{\alpha}:=	\beta_2 (\phi_\alpha^*(T)- \phi_\Omega) - \lambda\rho_\alpha .
\end{align}
Observe that the last summand of $f_\alpha$ is independent of time. By virtue of \eqref{boundal1}, {\bf (A6)}, and {\bf (A7)}, we have,
for every $\alpha\in (0,1]$,
\begin{align}
\label{boundfal}
&\|f_\alpha\|_{L^2(0,T;H)}\,+\,\|g_{\alpha}\|_{L^2(0,T;H)}
+ \|\rho_\alpha\|_{V}\,+\,\|\pi_{\alpha}\|_{V}\nonumber
\\\quad{}&\le\,
C\left(\|\phi^*_\alpha\|_{C^0([0,T];V)} + 
\|w^*_\alpha\|_{C^1([0,T];V)}\,+1\right)\,\le\,C,
\end{align}
where in the following we denote by $C$ positive constants that may depend on the data of the system but not on $\alpha\in (0,1]$.
 
According to \cite[Thm.~4.5]{CGSS4}, the adjoint system has under the assumptions \eqref{Uad}, \eqref{uminmax} and {\bf (A1)}--{\bf (A7)}
 a unique weak solution with the regularity
\begin{align}
	\label{regpal}
	\pal&\in H^1(0,T;V^*)\cap L^\infty(0,T;V)\cap L^2(0,T;W),
	\\
	\label{regqal}
	\qal & \in L^2(0,T;V),
	\\
	\label{regral}
	\ral & \in H^1(0,T;H)\cap L^\infty(0,T;V).
\end{align}
Moreover, by virtue of \cite[Thm.~4.7]{CGSS4}, we know that the first-order optimality condition for the optimal control $u^*_\alpha$ is characterized by the variational inequality
\begin{align}
\label{vugal}
\int_Q \left(\ral+\nu u^*_\alpha +(u^*_\alpha-u^*)\right)(v-u^*_\alpha)\,\ge\,0\quad\forall\,v\in\Uad\,.
\end{align}

The next step consists in passing to the limit as $\alpha\searrow0$ in both  the adjoint system \eqref{adsys1}--\eqref{adsys5}
and the variational inequality \eqref{vugal}. To this end, uniform (with respect to $\alpha\in (0,1]$) estimates for the
adjoint variables $(\pal,\qal,\ral)$ must be shown. A closer look at the system \eqref{adsys1}--\eqref{adsys5} reveals that 
there is an inherent difficulty. 
To this end, observe that \eqref{adsys2} and \eqref{adsys4} imply that ${\overline{\qal}(t)}=0$ for almost every $t\in (0,T)$.
 Therefore, testing of \eqref{adsys1} 
with the constant function $v=\boldsymbol{1}/|\Omega|$, integration with respect to time over $[t,T]$, and application of 
the well-known integration-by-parts rule for functions in $H^1(0,T;V^*)\cap L^2(0,T;V)$, using the endpoint 
conditions~\eqref{adsys5} along with the abbreviations \eqref{g1al}--\eqref{pial},
 yield the identity
\begin{align}
{\overline{\pal}(t)}\,&=
\overline{\pi_\alpha} + \lambda \overline{\rho_\alpha} 
\,-\,\lambda\,{\overline{\ral}(t)}\,+
\,\frac1{|\Omega|} \int_{Q^t} g_\alpha \nonumber\\
\label{meanpal}
&\quad\, -\,\frac 1{|\Omega|} \int_{Q^t} \bigl[ \gamma\pal\,+\,h''_\alpha(\phi^*_\alpha)\qal\,+\,F''(\phi^*_\alpha)\qal \bigr]\,.
\end{align}
Here, and for the remainder of this paper, we put
\begin{equation*}
Q^t:=\Omega \times (t,T) \quad\mbox{whenever $\,t\in [0,T)$.}
\end{equation*}
Apparently, the term $\,h''_\alpha(\phi^*_\alpha)\qal\,$ cannot be controlled. Indeed, although $\phi^*_\alpha$ satisfies the strict separation
condition \eqref{separation} for any fixed $\alpha>0$, a uniform bound cannot be expected, since it may well happen that the constants in
\eqref{separation} satisfy $r_*(\alpha)\searrow -1$ or $r^*(\alpha)\nearrow 1$ as $\alpha\searrow0$, in which case $h''_\alpha(\phi^*_\alpha)$
may become unbounded. Consequently, we cannot hope to pass to the limit as $\alpha\searrow0$ in the system \eqref{adsys1}--\eqref{adsys5} as it 
stands, not even in its weak form with test functions $v\in V$.
In order to overcome this difficulty, we employ an idea that goes back to \cite{CGS0}. To this end, recall that the mean value of $\qal$ vanishes
almost everywhere in $(0,T)$. Therefore, we deduce from \eqref{adsys2} and \eqref{adsys4} the identity 
\begin{equation}
\label{splitting}
\pal(t)-{\overline{\pal}(t)}={\cal N}\qal(t),
\end{equation}
with the operator ${\cal N}$ introduced in \eqref{defN}. Notice that $\pal\in H^1(0,T;V^*)$ and $\overline\pal\in H^1(0,T)$, whence we conclude 
that also ${\cal N}\qal\in H^1(0,T;V^*)$.  

The identity \eqref{splitting} enables us to eliminate $\overline\pal$ from the problem. For this purpose, we introduce the spaces
\begin{align}
\label{defH0V0}
H_0:=\{v\in H:\overline v=0\} \quad\mbox{and}\quad V_0:=\{v\in V:\overline v=0\}=V\cap H_0\,.
\end{align}
Now observe that the subspace $\mbox{span}\{ \boldsymbol{1}\}$ of spatially constant functions is the orthogonal complement of $H_0$ with respect
to the inner product of $H$; moreover, $H_0$ is a closed subspace of $H$ and therefore a Hilbert space itself when equipped with the standard
inner product in $H$. In addition, owing to the Poincar\'e--Wirtinger inequality~\eqref{poincare}, the expression
\begin{equation}
\label{equivnorm}
(v,w)_0:=\iO \nabla v\cdot \nabla w  \quad\mbox{for $v,w\in V_0$}
\end{equation}
defines an inner product on $V_0$ whose associated norm is equivalent to the standard norm of $V$. Obviously, $V_0$ becomes a Hilbert space when
endowed with the inner product $(\,\cdot\,,\,\cdot\,)_0$, and since, according to \cite[Cor.~5.3]{CGS0}, $V_0$ is densely embedded in $H_0$,
we can construct the Hilbert triple $(V_0,H_0,V_0^*)$ with the dense and compact embeddings $\,V_0\subset H_0 \subset V_0^*$ and the usual
identification that
\begin{equation}
\label{identify}
\langle v,w\rangle_{V_0} \,=\,(v,w) \quad\mbox{for all }\,v\in H_0 \,\mbox{ and }\,w\in V_0.
\end{equation}        
The idea now is to change the standard variational version of the system \eqref{adsys1}--\eqref{adsys5} by not admitting every $v\in V$ as 
test function in \eqref{adsys1}, but only those having zero mean value. In this way, 
we eliminate \,$\overline\pal$\, from the problem; indeed, we easily find that the pair 
$(\qal,\ral)$ solves the reduced system
\begin{align}
\label{adj1}                                          
&\langle -\dt {\cal N}\qal,v\rangle_{V_0} \,+\iO\nabla\qal\cdot\nabla v\,+\,\gamma\iO {\cal N}\qal\, v -\,\lambda\iO\dt\ral\,v\nonumber\\
&=\,-\iO \bigl(\hal''(\phi^*_\alpha)\qal +F''(\phi^*_\alpha)\qal \bigr)\,v\,+\iO g_{\alpha}\,v\nonumber\\
&\qquad\mbox{for all $v\in V_0$ and a.e. in $(0,T)$},
\\[1mm]
\label{adj2}
& -\iO\dt\ral \,v \,+\iO \nabla(\kappa_1\ral-\kappa_2(\boldsymbol{1}\circledast \ral))\cdot\nabla v\,-\,b\iO\qal\,v
\,=\iO f_\alpha v 
\nonumber\\
&\qquad\mbox{for all $v\in V$ and a.e. in $(0,T)$}, 
\\[1mm]
\label{adj3}
&{\cal N}\qal(T)\,=\,\pi_\alpha - \overline{\pi_\alpha} 
\,, \quad \ral(T)=\,\rho_\alpha.
\end{align}
At this point, it is worth observing that, because of the zero mean value condition, the space $V_0$ does not contain the space
$C^\infty_0(\Omega)$, and therefore the variational equality with test functions $v\in V_0$ cannot be interpreted as an equation in the
sense of distributions.

In the following, we deduce some a priori estimates for the reduced system \eqref{adj1}--\eqref{adj3}. Here we argue formally, where we
note that all of the following calculations can be performed rigorously on the level of an approximating Faedo--Galerkin system using as
basis functions the eigenfunctions $e_j$, normalized by $\|e_j\|=1$, for $j\in\enne$, of the Laplace operator 
with homogeneous Neumann conditions. That
is, we have 
\begin{align*}
&-\Delta e_j=\lambda_j e_j \quad\mbox{in }\,\Omega, \qquad \dn e_j=0 \quad\mbox{on }\,\Gamma, \qquad\mbox{for all\, $j\in\enne$},\\
&0=\lambda_1<\lambda_2\le\ldots,\quad \lim_{j\to\infty}\lambda_j=+\infty,\qquad {( e_i ,e_j)}=\delta_{ij}\quad\mbox{for all\, } i,j\in\enne.
\end{align*}
In this connection, observe that $\{e_j\}_{j\in\enne}$ forms a complete orthonormal system in $H$, while $\{e_j\}_{j\ge 2}$ is obviously
a complete orthonormal system in the space $H_0$ of functions having zero mean value, and the eigenspace associated with the eigenvalue
$\lambda_1=0$ is just the space of constant functions.

\vspace{2mm}\noindent
\underline{{\sc First estimate.}}
\quad 
We now insert $v=\qal(t)$ (which belongs to $V_0$) in \eqref{adj1}, and $ v=-\frac{\lambda}b \,\dt\ral(t)$ in \eqref{adj2} (this is only formal), and add the
resulting equations, whence a cancellation of two terms occurs. Then we integrate with respect to time over $(t,T)$ for arbitrary $t\in [0,T)$, 
taking \eqref{adj3} into account. 
Using the properties \eqref{dadefN}--\eqref{propN} of the operator ${\cal N}$, and rearranging terms, we obtain the identity
\begin{align}
\label{Juppuno}
&\frac 12\,\|\qal(t)\|^2_* \,+\int_{Q^t}|\nabla\qal|^2\,+\gamma\int_t^T \|\qal(s)\|^2_*\,{\rm d}s\,+\int_{Q^t} h''_\alpha(\phi_\alpha^*)\,|\qal|^2
\nonumber\\
&+\,\frac{\lambda}b\,\int_{Q^t}|\dt \ral|^2\,+\,\frac{\lambda\kappa_1}{2b}\,\|\nabla\ral(t)\|^2 - \frac 12\,\|\nabla \pi_\alpha \|^2 - \frac{\lambda\kappa_1}{2b}\,\|\nabla\rho_\alpha\|^2
\nonumber\\                                                                               
&= \,-\int_{Q^t} F''(\phi_\alpha^*)\,|\qal|^2\,+\int_{Q^t}g_\alpha\,\qal \,-\,\frac{\lambda}b \int_{Q^t} f_\alpha\,\dt\ral
\nonumber\\
&\quad\,\,+\,\frac{\lambda\kappa_2}{b}\int_{Q^t} \nabla(\boldsymbol{1}\circledast \ral)\cdot\dt\nabla\ral
\,=:\,I_1+I_2+I_3+I_4\,,
\end{align} 
with natural meaning. Observe that the fourth summand on the left-hand side is nonnegative. 
The last two terms on the left-hand side are uniformly bounded due to \eqref{boundfal}. 
Moreover, by virtue of \eqref{separation} and {\bf (A2)},
we have that \,\,$\|F''(\phi_\alpha^*)\|_{L^\infty(Q)}\,\le\,C$, and therefore it follows from Young's inequality, using \eqref{boundfal} and the compactness inequality
\eqref{compact}, that
\begin{align}
\label{Juppdue}
I_1+I_2\,\le\,C\,+\,C\int_{Q^t}|\qal|^2\,\le\,\frac 12\int_{Q^t}|\nabla\qal|^2\,+\,C\int_t^T\|\qal(s)\|^2_*\,{\rm d}s\,+\,C\,.
\end{align}
Moreover, by Young's inequality and \eqref{boundfal},
\begin{equation}
\label{Jupptre}
I_3\,\le\,\frac{\lambda}{2b}\int_{Q^t}|\dt\ral|^2\,+\,C\,.
\end{equation}
Finally, integration by parts with respect to time and Young's inequality yield the estimate
\begin{align}
\label{Juppquattro}
I_4\,&=\,-\,\frac{\lambda\kappa_2}b \iO \nabla(\boldsymbol{1}\circledast \ral)(t)\cdot\nabla\ral(t) \,+\,\frac{\lambda\kappa_2}b
\int_{Q^t}|\nabla\ral|^2
\nonumber\\
&\le\,\frac{\lambda\kappa_1}{2b}\,\|\nabla\ral(t)\|^2\,+\,C\,\|\nabla(\boldsymbol{1}\circledast\ral)(t)\|^2
\,+\,\frac{\lambda\kappa_2}b
\int_{Q^t}|\nabla\ral|^2\nonumber\\
&\le\,\frac{\lambda\kappa_1}{2b}\,\|\nabla\ral(t)\|^2\,+\,C\int_t^T\|\nabla\ral(s)\|^2\,{\rm d}s\,.
\end{align} 
Combining \eqref{Juppuno}--\eqref{Juppquattro}, and applying Gronwall's lemma backwards in time, we have thus shown the estimate
\begin{align}
\label{boundal2}
&\|\qal\|_{L^\infty(0,T;V^*)\cap L^2(0,T;V)}\,+\,\|\ral\|_{H^1(0,T;H)\cap L^\infty(0,T;V)}\nonumber\\
&+\int_{Q}h_\alpha''(\phi_\alpha^*)\,|\qal|^2\,\le\,C\,\quad \forall\,\alpha\in (0,1],
\end{align}
whence it obviously follows that
\begin{equation}
\label{boundal3}
\|\boldsymbol{1}\circledast \ral\|_{H^2(0,T;H)\cap W^{1,\infty}(0,T;V)}\,\le\,C \,\quad\forall\,\alpha\in (0,1].
\end{equation}
In addition, since it is known that the mapping ${\cal N}$ is a topological isomorphism between $V_0^*$ and $V_0$ and, for any $s\ge 0$,
between $H^s(\Omega)\cap H_0$ and $H^{s+2}(\Omega)\cap H_0$, we also have ${\cal N}\qal\in L^\infty(0,T;V_0)\cap L^2(0,T;H^3(\Omega))$ and  
\begin{equation}
\label{boundal4}
\|{\cal N}\qal\|_{L^\infty(0,T;V)\cap L^2(0,T;H^3(\Omega))}\,\le\,C\,\quad\forall\,\alpha\in (0,1].
\end{equation}

\vspace{2mm}\noindent
\underline{{\sc Second estimate.}} \quad
As a preparation for the next estimate, we introduce the space
\begin{equation}
\label{defZ}
{\cal Z}:=\{v\in H^1(0,T;V_0^*)\cap L^2(0,T;V_0):\,\,v(0)=0\},
\end{equation}
which, as a closed subspace, becomes a Hilbert space itself when endowed with the standard inner product
and norm of $H^1(0,T;V_0^*)\cap L^2(0,T;V_0)$. Notice that ${\cal Z} \subset C^0([0,T];H_0)$ which makes the 
initial condition $v(0)=0$ meaningful; we also have the dense and compact embeddings
$$\,{\cal Z}\subset L^2(0,T;H_0)\subset {\cal Z}^*.$$
Moreover, ${\cal Z}$ is dense in $L^2(0,T;V_0)$ since it 
contains the dense subspace $H^1_0(0,T;V_0)$. Therefore, the dual space $\,L^2(0,T;V_0^*){\cong}(L^2(0,T;V_0))^*$  
can be identified with a subspace of the dual space ${\cal Z}^*$ in the usual way, i.e., such that
\begin{equation}
\label{idenZ}
\langle v,w\rangle_{\cal Z}\,= \int_0^T\!\!\langle v(t),w(t)\rangle_{V_0}\,{\rm d}t
\quad \makebox{for all }\,v\in L^2(0,T;V_0^*)\,\mbox{ and }\,
w\in {\cal Z}\,.
\end{equation}
Now, we take an arbitrary $v\in {\cal Z}$ as test function in \eqref{adj1} and integrate over $(0,T)$. We obtain
\begin{align}
\label{pier1}                                          
&\int_0^T \langle -\dt {\cal N}\qal(t),v(t)\rangle_{V_0} {\rm d}t\,+\int_Q\nabla\qal\cdot\nabla v\,+\,
\int_Q \hal''(\phi^*_\alpha)\qal \, v
\nonumber \\
& =\int_Q \bigl(- \gamma \, {\cal N}\qal\,  + \,\lambda\, \dt\ral - F''(\phi^*_\alpha)\qal +g_{\alpha}\bigr)\,v\,.
\end{align}
{Next, we observe that 
${\cal N}\qal \in H^1(0,T;V_0^*)\cap L^2(0,T;V_0) 
%\quad \hbox{for all }\, \alpha\in (0,1]
$ 
and integrate by parts in the first term of \eqref{pier1}. With the help of \eqref{boundfal}, \eqref{adj3}, and \eqref{boundal4}, we infer that 
\begin{align}
\label{pier2}
&\int_0^T \!\langle - \dt {\cal N}\qal(t),v(t)\rangle_{V_0}\,{\rm d}t
%= |- \langle \dt {\cal N}\qal, v\rangle_{\cal Z}
= 
- (\pi_\alpha, v(T))  + \int_0^T\! \langle \dt v(t),{\cal N}\qal(t)\rangle_{V_0}\,{\rm d}t, 
\end{align}
and, consequently, for every $v\in {\cal Z}$ it holds
\begin{align}
&\Big |\int_0^T \langle - \dt {\cal N}\qal(t),v(t)\rangle_{V_0}\,{\rm d}t\Big |
\nonumber \\
&\le \, \| \pi_\alpha \|\, \|v\|_{C^0 ([0,T];H_0)} +
\int_0^T\|{\cal N}\qal(t)\|_{V_0}\,\|\dt v(t)\|_{V_0^*}\,{\rm d}t
\nonumber\\
&\le C\,\|v\|_{\cal Z} + C\,\|{\cal N}\qal\|_{L^2(0,T;V)}\,\|\dt v\|_{L^2(0,T;V_0^*)}
\le\,C\,\|v\|_{\cal Z}\,.
\label{pier3}
\end{align}
Hence, in view of the estimates \eqref{boundal2} and \eqref{pier3},
 we easily find from a comparison of terms in \eqref{pier1} that the linear functional
$$
\Lambda_\alpha:{\cal Z}\to\erre, \quad\langle\Lambda_\alpha,v\rangle_{\cal Z}:=\int_Q h_\alpha''(\phi^*_\alpha)\,\qal\,v\quad
\mbox{for }\,v\in {\cal Z},
$$
satisfies
\begin{equation}
\label{boundal6}
\|\Lambda_\alpha\|_{\cal Z^*}\,\le\,C \quad\forall\,\alpha\in (0,1].
\end{equation}

By the estimates shown above, there exist a sequence $\{\alpha_n\}_{n\in\enne}\subset (0,1]$ and limit points $q^*,r^*,\Lambda^*$ such that
$\alpha_n\searrow0$ and
\begin{align}
\label{limq}
q_{\alpha_n}\,&\,\to \,q^*\,&&\mbox{weakly star in }\, L^\infty(0,T;V^*)\cap L^2(0,T;V) ,
\\
\label{limr}
r_{\alpha_n}\,&\,\to \,r^*\,&&\mbox{weakly star in }\,H^1(0,T;H)\cap L^\infty(0,T;V),
\\
\label{lim1cr}
\boldsymbol{1}\circledast r_{\alpha_n}\,&\,\to\, \boldsymbol{1}\circledast r^*\,&&\mbox{weakly star in }\, H^2(0,T;H)\cap
W^{1,\infty}(0,T;V),
\\
\label{limNq}
{\cal N}q_{\alpha_n}\,&\,\to \,{\cal N}q^*\,&&\mbox{weakly star in }\, L^\infty(0,T;V_0)\cap L^2(0,T;H^3(\Omega)),
\\
\label{limlam}
\Lambda_{\alpha_n}\,&\,\to \,\Lambda^*\,&&\mbox{weakly in }\,{\cal Z}^*.
\end{align}
Moreover, in view of Theorem {\ref{THM:anktozero}} and Theorem {\ref{THM:APPROX}}, we may without loss of generality assume that
$u^*_{\alpha_n}\to u^*$ strongly in $L^2(Q)$ and that the convergence properties \eqref{conphial} and \eqref{conwal} for the state 
components $\phi^*_{\alpha_n}$ and $w^*_{\alpha_n}$ are satisfied correspondingly with $(\phi^0,w^0)=(\phi^*,w^*)$.
 Consequently, we have for $n\to\infty$ that
$$F''(\phi^*_{\alpha_n})q_{\alpha_n}\,\to\,F''(\phi^*)q^* \quad\mbox{weakly in }\,L^2(Q)$$ 
and that
\begin{align*}
&f_{\alpha_n}\to f^*:=\beta_3 ( \boldsymbol{1} \circledast (w^* - w_Q) )+\beta_4(w^*(T)-w_\Omega)
	+ \beta_5 (\dt w^* - w_Q'),\\[1mm]
&g_{\alpha_n}\to g^*:=\beta_1(\phi^*-\phi_Q),
\\[1mm]
&\rho_{\alpha_n} \to \rho^* := \beta_6 (\dt w^*(T)- w_\Omega') ,
\\[1mm]
&\pi_{\alpha_n} \to	\pi^* := \beta_2 (\phi^*(T)- \phi_\Omega) - \lambda\rho^* ,
\end{align*}
for suitable convergence properties as from \eqref{conphial} and \eqref{conwal}. 

Now we consider the system \eqref{adj1}--\eqref{adj3} for $\alpha=\alpha_n$, where we
replace \eqref{adj1} (and the first final condition in \eqref{adj3}) with the time-integrated
version \eqref{pier1}--\eqref{pier2}, with test functions $v\in {\cal Z}$. Passage to the limit as $n\to\infty$, using the above convergence
properties, then yields that
\begin{align}
\label{etze1}
&\langle \Lambda^*,v\rangle_{\cal Z}\,=\,
(\pi^*, v(T))  - \int_0^T \langle \dt v(t),{\cal N}q^*(t)\rangle_{V_0}\,{\rm d}t -\int_Q \nabla q^*\cdot\nabla v
\nonumber\\
&\quad{}+\int_Q\Bigl[{}-\gamma{\cal N}q^*+\lambda\dt r^*-F''(\phi^*)q^*+g^*\Bigr]v   \quad\mbox{for all $\,v\in{\cal Z}$},
\end{align}
\begin{align}
%\\[1mm]
\label{etze2}
& -\iO\dt r^*\,v \,+\iO \nabla(\kappa_1 r^*-\kappa_2(\boldsymbol{1}\circledast r^*))\cdot\nabla v\,-\,b\iO q^*\,v
\,=\iO f^*v 
\nonumber\\
&\qquad\mbox{for all $v\in V$ and a.e. in $(0,T)$}, 
\\[1mm]
\label{etze3}
&r^*(T)=\rho^*\,.
\end{align}

Finally, we consider the variational inequality \eqref{vugal} for $\alpha=\alpha_n$, $n\in\enne$. By passing to the
limit as $n\to\infty$, we find that
\begin{align}\label{limvug}
\int_Q (r^*+\nu \,u^*)(v-u^*)\,\ge\,0 \quad\forall\,v\in \Uad.
\end{align}

Summarizing the above considerations, we have proved the following first-order necessary optimality conditions for the
optimal control problem~{\bf (P)}.

\begin{theorem} 
Suppose that the conditions {\bf (A1)}--{\bf (A7)}, \eqref{Uad} and \eqref{uminmax} are fulfilled, and let $u^*\in\Uad$ be a 
minimizer of the optimal control problem {\bf (P)} with the associated uniquely determined state components \,$\phi^*,w^*$. 
Then there exist $\,q^*,r^*$, and $\Lambda^*$
such that the following holds true:
\begin{description}
\item[{\rm (i)}] \,\,\,$q^*\in L^\infty(0,T;V_0^*)\cap L^2(0,T;V)$, $\,{\cal N}q^*\in L^\infty(0,T;V_0)\cap L^2(0,T;H^3(\Omega))$, \hfill\break
%$\,\dt{\cal N}q^*\in {\cal Z}^*$,  
$r^*\in H^1(0,T;H)\cap L^\infty(0,T;V)$, $\,\boldsymbol{1}\circledast r^*\in
H^2(0,T;H)\cap W^{1,\infty}(0,T;V)$, 
$\Lambda^*\in {\cal Z}^*$.
\item[{\rm (ii)}] The adjoint system \eqref{etze1}--\eqref{etze3} and the variational inequality \eqref{limvug}
are satisfied.
\end{description} 
\end{theorem}

\vspace{2mm}
\begin{remark}{\em
{(i)} Observe that the adjoint state variables $(q^*,r^*)$ and the Lagrange multiplier $\Lambda^*$ are not uniquely
determined. However, all possible
choices satisfy \eqref{limvug}, i.e., $u^*\,$ is for $\nu>0$ the $L^2(Q)$-orthogonal projection of $-{\nu}^{-1} r^*$ onto the 
closed and convex set $\Uad$, and for a.e. $\,(x,t) \in Q$ it holds 
\begin{align}
	\non
	u^*(x,t)=\max \big\{ 
	u_{\rm min}(x,t), \min\{u_{\rm max}(x,t),-{\nu}^{-1} r^*(x,t)\} 
	\big\}. 
\end{align} 
{(ii)} We have, for every $n\in\enne$, the complementarity slackness condition}
$$\langle \Lambda_{\alpha_n},q_{\alpha_n}\rangle_{{\cal Z}}=\int_Q h_{\alpha_n}''(\phi^*_{\alpha_n})\,|q_{\alpha_n}|^2
=\int_Q \frac{2\alpha_n}{1 - (\phi^*_{\alpha_n})^2}\,|q_{\alpha_n}|^2
\,\ge\,0.
$$
{\em Unfortunately, our convergence properties for  $\{\phi^*_{\alpha_n}\}$
and $\{q_{\alpha_n}\}$ do not permit a passage to the limit in this
inequality to
derive a corresponding result for {\bf (P)}.}
\end{remark}

%%%%%%%%%%%%%%%%%%%%%%%%%%%%%%%%%%%%%%%%%%%%%%%%%%%%%%%%%%%%%%%%%%%%%%%%%%%

\vskip 3mm
\noindent{\bf Acknowledgements}

\noindent
This research was supported by the Italian Ministry of Education, 
University and Research (MIUR): Dipartimenti di Eccellenza Program (2018--2022) 
-- Dept.~of Mathematics ``F.~Casorati'', University of Pavia. 
In addition, {PC and AS gratefully acknowledge some other support 
from the MIUR-PRIN Grant 2020F3NCPX ``Mathematics for industry 4.0 (Math4I4)'' and}
their affiliation to the GNAMPA (Gruppo Nazionale per l'Analisi Matematica, 
la Probabilit\`a e le loro Applicazioni) of INdAM (Isti\-tuto 
Nazionale di Alta Matematica). 

%%%%%%%%%%%%%%%%%%%%%%%%%%%%%%%%%%%%%%%%%%%%%%%%%%%%%%%%%%%%%%%%%%%%%%%%%%%
\enddocument